\documentclass[12pt]{amsart} 

\usepackage{amssymb} 
\usepackage{amsmath} 
\usepackage{fullpage} 
\usepackage{url}

\newtheorem{theo}{Theorem}
\newtheorem{prop}{Proposition}
\newtheorem{lemma}{Lemma}

\theoremstyle{definition}

\theoremstyle{remark}
\newtheorem{rem}{Remark}

\def\Er{{\mathbb E}}
\def\Fr{{\mathbb F}}

\def\Nr{{\mathbb N}}

\def\Pr{{\mathbb P}}
\def\Qr{{\mathbb Q}}
\def\Rr{{\mathbb R}}

\def\Ac{{\mathcal{A}}}
\def\Bc{{\mathcal{B}}}

\def\Fc{{\mathcal{F}}}

\def\Sc{{\mathcal{S}}}

\def\one{{\rm \bf 1}}
\def\({\left(}     
\def\){\right)}    
\def\[{\left[}     
\def\]{\right]}

\def\as{{\frenchspacing a.s.}~}
\def\ep{\varepsilon} 
\begin{document}
\title{Law of Large Numbers for Risk Measures  \footnote{{\it AMS-}Classification 90B50, 91B06, 91B16, 91G99} } 
\author{Freddy Delbaen}
\address{Departement f\"ur Mathematik, ETH Z\"urich, R\"{a}mistrasse
   101, 8092 Z\"{u}rich, Switzerland}
 \address{Institut f\"ur Mathematik,
 Universit\"at Z\"urich, Winterthurerstrasse 190,
 8057 Z\"urich, Switzerland}
\date{First version April 26, 2020, this version \today}

\begin{abstract}
Under appropriate integrability conditions the risk measure of the sample measures for a law  invariant risk measure converge almost surely to the risk measure of the sampled random variable. The results follow from general convergence theorems based on the theory of Orlicz spaces.\end{abstract}

\maketitle

\section{Introduction and Notation}
The notation we use is the standard notation of probability theory.  $(\Omega,\Fc,\Pr)$ will denote a probability space without atoms.  If $\Qr$ is a probability defined on $\Fc$ and absolutely continuous with respect to $\Pr$, we will identify $\Qr$ with its Radon-Nikod\'ym derivative $h=d\Qr/d\Pr$. A basic example of such an atomless space is $([0,1],\Bc,m)$ where $\Bc$ is the Borel sigma algebra and $m$ is Lebesgue measure.  This space will also play a fundamental role in the sequel.

Because of neganumerophobia we will use a setup closer to premium calculation.  For a convex closed set of probability densities, $\Sc\subset L^1$ we denote by $\rho$ the functional (whenever defined)
$$
\rho(\xi)=\sup\{\Er[\xi h\mid h\in \Sc\}.
$$
This functional is up to sign changes the same as a coherent monetary utility function.  The interpretation of $\rho$ is the same as for the utility functions. Among the many uses we mention  risk measures and insurance premium principles defined on random variables. A basic problem is the following.  Suppose we see a large sample of a random variable.  Can we define the functional $\rho$ on such samples and does this procedure converge to the risk measure of the random variable?  The best is to put the problem in the framework of risk measures that are law invariant (law determined or rearrangement invariant).  That means that when $\eta$ and $\eta'$ have the same distribution, then $\rho(\eta)=\rho(\eta')$.  Since we will be working on atomless probability spaces, this property is also equivalent to the property that $\Sc$ is rearrangement invariant.  Law determined risk measures can be defined for probability distributions on $\Rr$ with compact support (the latter to avoid integrability problems).  It therefore does not matter on which probability space we work. In later sections we will start with a risk measure defined for random variables defined on $\Omega$ but then we will ``transport" it to random variables defined on $[0,1]$. For more information on rearrangement invariant measures see \cite{kusuoka} where the general form of these risk measures is given.

The risk measures we will use are sometimes defined on bigger spaces than $L^\infty$, the space of bounded random variables.  The best way to study this is to use  Orlicz spaces.  For the theory of these spaces we refer to \cite{Rao-Ren}.  In this book the reader will find a very general approach.  We will restrict ourselves to Young functions that are more regular than the ones in \cite{Rao-Ren}.  For our presentation this restriction is of no importance but it avoids a  lot of -- in our opinion cosmetic --problems. When we say $\Phi$ is a Young function, we mean that $\Phi$ is an even real valued function defined on $\Rr$ with the help of its derivative $\phi$.  This derivative function $\phi$ will be nondecreasing on $\Rr_+$, continuous at $0$, $\phi(x)=0$ if and only if $x=0$ and $\lim_{x\rightarrow +\infty}\phi(x)=+\infty$.  We do not require more regularity for $\phi$. The conjugate Young function $\Psi$ is then defined using the right continuous inverse function, $\psi$, of $\phi$. More precisely for $x\ge 0$: $\psi(x)=\inf\{u\mid \phi(u)> x\}$. Hence $\Psi(x)=\int_0^{|x|}\psi(u)\,du$.  We realise that the class of Orlicz spaces so defined, does not contain $L^1$ or $L^\infty$, so when appropriate we will make additional remarks.  We follow \cite{Rao-Ren} and define
\begin{eqnarray*}
L^\Phi&=&\left\{ \xi \mid \text{ there is $k>0$ with } \Er\[ \Phi(k\xi)\]<\infty  \right\}\\
M^\Phi&=&\left\{ \xi \mid \text{ for all $k>0$ we have } \Er\[ \Phi(k\xi)\]<\infty  \right\}=L^{(\Phi)}\subset L^\Phi.
\end{eqnarray*}
These spaces are usually equipped with the Luxemburg norm
$$
\Vert \xi \Vert_\Phi=\inf \left\{  \alpha >0 \mid \Er\[\Phi\(\frac{\xi}{\alpha}\)  \]\le 1 \right\}.
$$
The dual space of $M^\Phi$ is the space $L^\Psi$ and we have for $\xi\in L^\Phi,\eta\in L^\Psi$ that $\xi\eta\in L^1$ and $|\Er[\xi\eta]|\le 2 \Vert \xi\Vert_\Phi \Vert\eta\Vert_\Psi$. The dual space of $L^\Phi$ is not necessarily $L^\Psi$.  Indeed we have that $(L^\Phi)^*=L^\Psi$ if and only if $L^\Phi=M^\Phi$ if and only if $\Phi$ satisfies the $\Delta_2$ condition:  there is $x_0$ and a constant $C$ such that for all $x\ge x_0$ the inequality $\Phi(2x)\le C \Phi(x)$ is valid.  If this is the case one can modify the function $\Phi$ without changing the spaces and such that in the inequality we may suppose that $x_0=0$. 
\section{Main Result}
In this section we suppose that the scenario set $\Sc$, for the  risk measure $\rho$,  is supposed to be a subset of $L^\Phi$. The set $\Sc$ is a convex closed subset of $L^1$ and hence it must be bounded in $L^\Phi$.  Indeed, if $h_n\in \Sc$ is such that $\Vert h_n\Vert_\Phi$ is an unbounded sequence, say $\Vert h_n\Vert_\Phi\ge 4^n$, then the element $h=\sum_n \frac{h_n}{2^n}$ is still in $\Sc$ but we have for each $n$: $\Vert h\Vert_\Phi\ge \Vert h_n 2^{-n}\Vert_\Phi \ge 2^{-n} 4^n=2^n$.  This would mean that $\Sc$ is not a subset of $L^\Phi$.   The set $\Sc$ is also weak$^*$, meaning $\sigma(L^\Phi,M^\Psi)$, compact. To see this we must only show that it is weak$^*$ closed.  This follows from the (obvious) continuity of the injection $L^\Phi\rightarrow L^1$ for the topologies $\sigma(L^\Phi,M^\Psi)\text{ and } \sigma(L^1,L^\infty)$, just observe that $L^\infty\subset M^\Phi$.  This continuity also implies that $\Sc$ is closed in $L^\Phi$. We remark that when $\Sc\subset L^1$ is a convex set of probabilities that is closed in $L^\Phi$, this does not necessarily imply that $\Sc$ must be closed in $L^1$.

Because of this boundedness, the functional $\rho$ with scenario set $\Sc$ can be defined on $L^\Psi$.  Indeed for $\eta\in L^\Psi$ we have 
$$
| \rho(\eta)|=|\sup\{\Er[\eta h]\mid h\in \Sc\}| \le  2\Vert \eta\Vert_\Psi \sup_{h\in\Sc} \Vert h\Vert_\Phi < \infty.
$$
As expected this boundedness also implies that $\rho$ is continuous on $L^\Psi$.  For completeness we include the obvious proof.
\begin{lemma} If $\eta_n$ is a sequence in $L^\Psi$ tending to $\eta$ for $\Vert . \Vert_\Psi$, then $\rho(\eta_n)\rightarrow \rho(\eta)$. More precisely $|\rho(\eta)-\rho(\eta')|\le 2 \sup_{h\in\Sc}\Vert h\Vert_\Phi \Vert \eta-\eta'\Vert_\Psi $.
\end{lemma}
{\bf Proof. }  We only have to prove the inequality.  For given $\ep>0$ and $\eta$ we choose $h\in \Sc$ with $\rho(\eta)\le \Er[h\eta]+\ep$.  This gives
$$
\rho(\eta)\le \Er[h\eta]+\ep = \Er[h\eta']+\Er[h(\eta-\eta')]+\ep\le\rho(\eta')+2\Vert h\Vert_\Phi\Vert \eta-\eta'\Vert_\Psi+\ep.
$$
Because $\ep>0$ was arbitrary we get $\rho(\eta)\le \rho(\eta')+2 \sup_{h\in\Sc}\Vert h\Vert_\Phi \Vert \eta-\eta'\Vert_\Psi $.  We can interchange $\eta$ and $\eta'$ and then get the desired result.\qed

\medskip
{\it From now on we suppose that $\rho$ is  rearrangement invariant. }
\medskip

For a probability measure $\mu$ on $\Rr$, we use the notation $\rho(\mu)$ for the outcome $\rho(\eta)$, where $\mu$ is the distribution of $\eta$.  When $\Sc\subset L^\Phi$, this is defined as soon as $\int_\Rr \Psi(x)d\mu<\infty$, i.e. $\eta\in L^\Psi$. 

For probability laws on $\Rr$ there is an easy way to transform convergence of probability distributions into convergence of random variables (a special case of the Skorohod lemma). Suppose that $\mu_n$ converges weak$^*$ to $\mu$, i.e. for each bounded continuous function $g$ defined on $\Rr$ we have $\int_\Rr g(x)d\mu_n\rightarrow \int_\Rr g(x)d\mu$.  The distribution functions are defined as $\Fr_n(x)=\mu_n(-\infty, x]$, $\Fr(x)=\mu(-\infty,x]$.  Weak$^*$ convergence is then equivalent to $\Fr_n(x)\rightarrow \Fr(x)$ for each continuity point $x$ of $\Fr$. The quantile functions are defined as $q_n(x)=\inf\{u\mid \Fr_n(u)>x\}$, $q(x)=\inf\{u\mid \Fr(u)>x\}$.  The quantile functions, $q_n$, are defined on the unit interval $(0,1)$ and under the Lebesgue measure they have the distribution functions $\mu_n$. The convergence of $\mu_n$ is then transformed into almost everywhere convergence of the quantile functions. Convergence $q_n\rightarrow q$ in the Orlicz space space $L^\Psi$ is not easily defined in terms of the measures $\mu_n$. The following is a partial answer to this question
\begin{prop} Suppose that $\mu_n$ converges weak$^*$ to $\mu$.  If for each $k>0$, $\int_\Rr\Psi(kx)\mu_n(dx)\rightarrow \int_\Rr\Psi(kx)\mu(dx)<\infty$, we have $\Vert q_n-q\Vert_\Psi\rightarrow 0$.  Consequently $\rho(\mu_n)\rightarrow \rho(\mu)$.
\end{prop}
{\bf Proof. }  We have $\rho(\mu_n)=\rho(q_n),\rho(\mu)=\rho(q)$ so the last line follows from the lemma. To prove the norm convergence we have to estimate $\int_{[0,1]}\Psi(k|q_n(u)-q(u)|)\,du$. We use the inequalities
$$
\int \Psi(k|q_n-q|)\le \int \Psi(2k\max(q_n,q))\le \int\Psi(2k q_n)+\int\Psi(2kq).
$$
Because $q_n\rightarrow q$ \as and $\int\Psi(2kq_n)\rightarrow \int\Psi(2kq)$, Scheff\'e lemma implies that there is $n_0$ such that $\Psi(2k q_n), n\ge n_0$ forms a uniformly integrable family.  Consequently the elements $\int \Psi(k|q_n-q|)$ (for $n\ge n_0$) form a uniformly integrable family.  Therefore for every $k$: $\int_{[0,1]}\Psi(k|q_n(u)-q(u)|)\,du\rightarrow 0$.  This shows that $\lim_n \Vert q_n-q\Vert_\Psi \rightarrow 0$.\qed
\begin{rem} The hypothesis implies that $q\in M^\Psi$ but there is no need to suppose that $q_n\in M^\Psi$.
\end{rem}

We are now ready to state the main result of this section.

\begin{theo} We use the notation and assumptions from the beginning of this section.  Let $\xi\in M^\Psi$.  Let $\xi_n\colon (\Omega,\Ac,\Pr)\rightarrow \Rr$ be an iid sequence with the same distribution $\mu$ as $\xi$.  For each $N$ we define the sample distribution as $\mu_n=\frac{1}{N}\sum_{n=1}^N\delta_{\xi_n}$ where $\delta_y$ denotes the Dirac measure concentrated in $y$.  Then almost surely $\rho(\mu_N)\rightarrow \rho(\mu)$.
\end{theo}
{\bf Proof. } The usual law of large numbers allows to conclude that almost everywhere the random measures $\frac{1}{N}\sum_{n=1}^N\delta_{\xi_n}$ converge weak$^*$ to $\mu$. The variant of Skorohod's theorem then permits to find for almost every $\omega \in \Omega$, increasing functions $q_N$ on $[0,1]$ such that $q_N$ has as its distribution $\mu_n=\frac{1}{N}\sum_{n=1}^N\delta_{\xi_n}$.  They converge almost everywhere (on $[0,1]$) to $q$,  the increasing rearrangement of $\xi$. Again the law of large numbers shows that almost surely for every $0<K\in\Nr$,  $\frac{1}{N}\sum_{n=1}^N {\Psi}(K \xi_n)$ converges to $\Er[{\Psi}(K\xi)]$ (where the set of measure 1 can be chosen to be independent of $K$). On $[0,1]$ this translates as $\int_{[0,1]}\Psi(Kq_N)\rightarrow \int_{[0,1]}\Psi(Kq)$.  We now apply the previous proposition.
 \qed

\begin{rem} With minor notational modifications, the above proof also works for the case $\Psi(x)=x$, i.e. the $L^1$ space.  In particular it also proves the theorem for the tail expectation.  In this case the result was proved by Van Zwet.
\end{rem}
\section{A more general Result}
In the section we generalise the result of the preceding section to the case of random variables in $L^\Psi$ (and not just $M^\Psi$).  We need an additional hypothesis on the scenarioset.  We recall the result of Ando, \cite{Ando}, which characterises the $\sigma(L^\Phi,L^\Psi)$ relative sequentially compact sets.
\begin{theo} (Ando) A set $A\subset L^\Phi$ is $\sigma(L^\Phi,L^\Psi)$ relative sequentially compact if and only if
$$
\lim_{\lambda\rightarrow \infty}\sup_{h\in A}\lambda\Er\[\Phi\(\frac{h}{\lambda}\)\]=0.
$$
In particular if $\Psi$ is $\Delta_2$ then bounded sets $A\subset L^\Phi$ automatically satisfy the above limit condition.
\end{theo}
We use the same notation as in the previous sections. In particular $\rho$ is defined on $L^\Psi$ and is rearrangement invariant.
\begin{prop} Suppose that $\Sc$ satisfies the conditions of Ando's theorem. Suppose that $\mu_n$ converges weak$^*$ to $\mu$.  If for some $k>0$, $\int_\Rr\Psi(kx)\mu_n(dx)\rightarrow \int_\Rr\Psi(kx)\mu(dx)<\infty$, we have $\sup_{h\in\Sc}\Er[h|q_n-q|]\rightarrow 0$.  Consequently $\rho(\mu_n)\rightarrow \rho(\mu)$.
\end{prop}
{\bf Proof. } By homogeneity we may suppose that $k=2$ and that 
$\int_\Rr\Psi(2x)\mu_n(dx)\rightarrow \int_\Rr\Psi(2x)\mu(dx)$. We may also suppose that for all $n$: $\int_\Rr\Psi(2x)\mu_n(dx)<\infty$. We replace the integrals on $\Rr$ by integrals of the quantile functions and get that $\int_{[0,1]}\Psi(2q_n)\rightarrow \int_{[0,1]}\Psi(2q)$.  Scheff\'e's lemma then implies that the sequence $\Psi(2q_n)$ is uniformly integrable. The obvious inequality $\Psi(|q_n-q|)\le \Psi(2q_n)+\Psi(2q)$ then implies that the sequence $\Psi(|q_n-q|)$ is also uniformly integrable. Hence $\int_{[0,1]}\Psi(|q_n-q|)\rightarrow 0$, since $q_n\rightarrow q$ \as.  We now have
$$
\sup_{h\in \Sc}\int_{[0,1]} h|q_N-q|\rightarrow 0.
$$ 
For given $\epsilon>0$ we take $\lambda$ so that $\lambda\int_{[0,1]}\Phi(\frac{h}{\lambda})\le \epsilon$.  Then we take $N_0$ so that for $N\ge N_0$: $\lambda\int_{[0,1]}\Psi(|q_N-q |)\le \epsilon$.  Young's inequality now yields $h|q_N-q|\le\lambda\Phi(\frac{h}{\lambda})+\lambda\Psi(|q_N-q |)$ and after integration we get for $N\ge N_0$:
$$
\sup_{h\in \Sc}\int_{[0,1]} h|q_N-\eta|\le 2\epsilon.
$$
\qed

\begin{theo} Suppose that the scenarioset $\Sc$ satisfies the condition of Ando's theorem. Let $\xi\in L^\Psi$.  Let $\xi_n\colon (\Omega,\Ac,\Pr)\rightarrow \Rr$ be an iid sequence with the same distribution $\mu$ as $\xi$.  For each $N$ we define the sample distribution as $\mu_n=\frac{1}{N}\sum_{n=1}^N\delta_{\xi_n}$ where $\delta_y$ denotes the Dirac measure concentrated in $y$.  Then almost surely $\rho(\mu_N)\rightarrow \rho(\mu)$.
\end{theo}
{\bf Proof.} We already know that almost surely $\mu_N$ converges weak$^*$ to $\mu$. Without loss of generality we may replace $\xi$ by a multiple of it and hence we may suppose that $\Er\[\Psi(2 \xi)\]\le 1$.  The law of large numbers implies that:
$$
\frac{1}{N}\sum_{k=1}^N \Psi(2\xi_k)\rightarrow \Er[\Psi(2\xi)],
$$
almost  surely and in $L^1$. We can now apply the previous proposition.\qed

\section{The case of Choquet Integration}
Choquet integration provides us with one of the most used examples of rearrangement invariant risk measures. As proved by David Schmeidler, \cite{Schm2}, \cite{FDbook}, a rearrangement invariant risk measure that is commonotone is given by the core of a convex game and hence uses a scenarioset of the form
$$
\Sc=\left\{h\in L^1\mid \Er[h]=1\text{ and for each }A\in\Fc:\Er[h\one_A]\ge f(\Pr[A])   \right\},
$$
where $f$ is a convex function $f\colon [0,1]\rightarrow [0,1],f(0)=0,f(1)=1$. In our context, where we want the risk measure to be defined on a larger space than $L^\infty$ we must use a weakly ($\sigma(L^1,L^\infty)$) compact set $\Sc$. This is equivalent to the left continuity of $f$ at the point $1$. Ryff's theorem, \cite{Ryff}, \cite{FDbook}, then says that $\Sc$ is the $L^1-$closed convex hull of random variables having the same law as the derivative $f'$ of $f$.  The statement $\Sc\subset L^\Phi$ is equivalent to $f'\in L^\Phi$. One can show that when $f'\in L^\Phi[0,1]$, the closed convex hull taken in $L^\Phi$ is the same as when the closedness is taken in $L^1$. The set $\Sc$ is also the set of those nonnegative random variables $h$ such that for each convex function $\beta\colon \Rr_+\rightarrow \Rr$ we have $\Er[\beta(h)]\le\int_{[0,1]}\beta(f')$ (convex dominance). Because of  this dominance the set $\Sc$ automatically satisfies Ando's criterium for $\sigma(L^\Phi,L^\Psi)$ sequential compactness.  As a result the topologies $\sigma(L^\Phi,L^\Psi)$ and $\sigma(L^1,L^\infty)$ coincide on the set $\Sc$. We get that the statistical results of the previous section apply to this case.  For iid sequences we find the almost sure convergence
$$
\xi_{[n:1]}f\(\frac{1}{n}\)+\sum_{k=1}^{k=n-1}\xi_{[n:k+1]}\(f\(\frac{k+1}{n}\)-f\(\frac{k}{n}\)\)=\int_0^1 q_n\,f'\rightarrow \int_0^1 q f'=\rho(\xi)
$$
where $\xi_{[n:1]},\ldots,\xi_{[n:n]}$ are the observations $\xi_1,\ldots,\xi_n$ in increasing order.  This result can also be obtained using  integration theory.

\end{document}